\documentclass[12pt,a4paper]{article}

\usepackage{amsmath}
\usepackage{amssymb}
\usepackage{enumerate}
\usepackage{t1enc}
\usepackage[latin1]{inputenc}
\usepackage[english]{babel}
\usepackage{amsfonts}
\usepackage{latexsym}
\usepackage{bm}

\newtheorem{theorem}{Theorem}[section]

\newtheorem{lemma}[theorem]{Lemma}

\newcommand\one{{\bf1}}

\DeclareMathOperator{\ourmod}{mod^*}

\newcommand*{\qed}{\hfill\ensuremath{\Box}}

\begin{document}
\title{Intersecting generalised permutations}

\author{Peter Borg\\[5mm]
Department of Mathematics \\
University of Malta\\
\texttt{peter.borg@um.edu.mt} 
\and
Karen Meagher \\ [5mm]
Department of Mathematics and Statistics \\
University of Regina \\
\texttt{karen.meagher@uregina.ca}
}

\date{} \maketitle

\begin{abstract}
  For any positive integers $k,r,n$ with $r \leq \min\{k,n\}$, let
  $\mathcal{P}_{k,r,n}$ be the family of all sets $\{(x_1,y_1), \dots,
  (x_r,y_r)\}$ such that $x_1, \dots, x_r$ are distinct elements of
  $[k] = \{1, \dots, k\}$ and $y_1, \dots, y_r$ are distinct elements
  of $[n]$. The families $\mathcal{P}_{n,n,n}$ and $\mathcal{P}_{n,r,n}$
  describe \emph{permutations} of $[n]$ and \emph{$r$-partial
    permutations} of $[n]$, respectively. If $k \leq n$, then
  $\mathcal{P}_{k,k,n}$ describes permutations of $k$-element subsets of
  $[n]$. A family $\mathcal{A}$ of sets is said to be
  \emph{intersecting} if every two members of $\mathcal{A}$
  intersect. In this note we use Katona's elegant cycle method to show
  that a number of important Erd\H{o}s-Ko-Rado-type results by various
  authors generalise as follows:~the size of any intersecting
  subfamily $\mathcal{A}$ of $\mathcal{P}_{k,r,n}$ is at most ${k-1
    \choose r-1}\frac{(n-1)!}{(n-r)!}$, and the bound is attained if
  and only if $\mathcal{A} = \{A \in \mathcal{P}_{k,r,n} \colon (a,b)
  \in A\}$ for some $a \in [k]$ and $b \in [n]$.
\end{abstract}

\section{Introduction}

For an integer $n \geq 1$, the set $\{1, 2, \dots, n\}$ is denoted by
$[n]$. The \emph{power set} $\{A \colon A \subseteq X\}$ of a set
$X$ is denoted by $2^X$, and the \emph{uniform} subfamily $\{Y
\subseteq X\colon |Y| = r\}$ of $2^X$ is denoted by $X \choose
r$. We call a set of size $n$ an \emph{$n$-set}.

If $\mathcal{F}$ is a family of sets and $x$ is an element in the
union of all sets in $\mathcal{F}$, then we call the subfamily of all
the sets in $\mathcal{F}$ that contain $x$ the \emph{star of
  $\mathcal{F}$} with \emph{centre} $x$. A family $\mathcal{A}$ is
said to be \emph{intersecting} if $A \cap B \neq \emptyset$ for every $A, B \in \mathcal{A}$. Note that a star of a family is intersecting.

The classical Erd\H{o}s-Ko-Rado (EKR) Theorem~\cite{EKR} says that
if $r \leq n/2$, then an intersecting subfamily $\mathcal{A}$ of
${[n] \choose r}$ has size at most ${n-1 \choose r-1}$, i.e. the
size of a star of ${[n] \choose r}$. If $r < n/2$, then, by the
Hilton-Milner Theorem \cite{HM}, $\mathcal{A}$ attains the bound
if and only if $\mathcal{A}$ is a star of ${[n] \choose r}$. Two
alternative proofs of the EKR Theorem that are particularly short
and beautiful were obtained by Katona \cite{K} and Daykin
\cite{D}. In his proof, Katona introduced a very elegant
technique called the \emph{cycle method}. Daykin's proof is based
on a fundamental result known as the Kruskal-Katona Theorem
\cite{K,Kr}. The EKR Theorem inspired a wealth of results and
continues to do so; see \cite{Borg7, DF,F}.

For positive integers $k,r,n$ with $r \leq \min\{k,n\}$, let
\begin{align*} 
\mathcal{P}_{k,r,n} := \{\{(x_1,y_1), \dots, (x_r,y_r)\} \colon &x_1, \dots, x_r \mbox{ are distinct elements of } [k], \\
& y_1, \dots, y_r \mbox{ are distinct elements of } [n]\}. 
\end{align*}
We shall call $\mathcal{P}_{k,r,n}$ a family of \emph{generalised
  permutations}. This is due to the fact that the elements of
$\mathcal{P}_{n,n,n}$ are permutations of the set $[n]$; the permutation
$y_1y_2 \dots y_n$ of $[n]$ corresponds uniquely to the set
$\{(1,y_1), (2,y_2), \dots, (n,y_n)\}$ in $\mathcal{P}_{n,n,n}$.

In the more general case where $k \leq n$, the set
$\mathcal{P}_{k,k,n}$ describes permutations of $k$-subsets of $[n]$;
a permutation $y_1y_2 \dots y_k$ of a $k$-subset of $[n]$ corresponds
uniquely to the set $\{(1,y_1), (2,y_2), \dots, (k,y_k)\}$ in
$\mathcal{P}_{k,k,n}$. The set $\mathcal{P}_{n,r,n}$ describes
\emph{$r$-partial permutations} of $[n]$ (see \cite{KL}).  The ordered
pairs formulation we are using follows \cite{Borg} and also
\cite{Borg1, Borg2}, in which very general frameworks are considered.

In the case $r = k$, if two sets $\{(1,y_1), (2,y_2), \dots,
(k,y_k)\}$ and $\{(1,z_1), (2,z_2),\dots, (k,z_k)\}$ in
$\mathcal{P}_{k,k,n}$ intersect, then $y_i = z_i$ for some $i \in
[k]$, and this is exactly what we mean by saying that the
permutations $y_1y_2 \dots y_k$ and $z_1z_2 \dots z_k$ (of two
$k$-subsets of $[n]$) intersect. In general, two generalised
permutations intersect if and only if they have at least one
ordered pair in common.

In this note we are concerned with the EKR problem for generalised
permutations.  We need only to consider the problem with $k \leq n$.
To see this, define $\lambda \colon [k] \times [n] \rightarrow [n]
\times [k]$ by $\lambda(x,y) := (y,x)$, then $\Lambda \colon
\mathcal{P}_{k,r,n} \rightarrow \mathcal{P}_{n,r,k}$ by
\[
\Lambda(\{(x_1,y_1),\dots, (x_r,y_r)\}) := 
  \{\lambda(x_1,y_1),\dots, \lambda(x_r,y_r)\} = \{(y_1,x_1),\dots,(y_r,x_r)\}.
\]
The functions $\lambda$ and $\Lambda$ are clearly both bijections.
Moreover, any $P,Q \in \mathcal{P}_{k,r,n}$ are intersecting if and
only if $\Lambda(P), \Lambda(Q)\in \mathcal{P}_{r,k,n}$ are
intersecting.  Therefore, throughout the rest of the paper it is to be
assumed that $k \leq n$.

The origins of our problem lie in \cite{DF1}, in which Deza and Frankl
prove that the size of an intersecting family of permutations of $[n]$
is at most $(n-1)!$, i.e. the size of a star of
$\mathcal{P}_{n,n,n}$. Cameron and Ku \cite{CK} extended this result
by establishing that only the stars of $\mathcal{P}_{n,n,n}$ attain
the bound (other proofs of this result are found in \cite{BH, GM, Li,
  WZ}).  This result was also done independently by Larose and
Malvenuto \cite{LM}, who actually showed that the stars of
$\mathcal{P}_{k,k,n}$ are the largest intersecting subfamilies of
$\mathcal{P}_{k,k,n}$ (see \cite[Theorem 5.1]{LM}). These results
summarize as follows.

\begin{theorem}[\cite{CK,DF1,LM}] \label{LM} The size of any
intersecting subfamily of $\mathcal{P}_{k,k,n}$ is at most
$\frac{(n-1)!}{(n-k)!}$, and the bound is attained only by the stars
of $\mathcal{P}_{k,k,n}$.
\end{theorem}

Ku and Leader \cite{KL} solved the EKR problem for $r$-partial
permutations of $[n]$ using Katona's cycle method. Moreover, they
showed that for $8 \leq r \leq n-3$, the largest intersecting
subfamilies of $\mathcal{P}_{n,r,n}$ are the stars. They
conjectured that only the stars are extremal for the few
remaining values of $r$ too. A proof of this conjecture, also
based on the cycle method, was obtained by Li and Wang \cite{LW}.

\begin{theorem}[\cite{KL,LW}]\label{KL} 
  For $r \in [n-1]$, the size of any intersecting subfamily of
  $\mathcal{P}_{n,r,n}$ is at most $ {n-1 \choose
    r-1}\frac{(n-1)!}{(n-r)!}$, and the bound is attained only by the
  stars of $\mathcal{P}_{n,r,n}$.
\end{theorem}

The scope of this note is to show that the methods used in
\cite{KL, LW} allow us to generalise Theorems~\ref{LM} and
\ref{KL} as follows.

\begin{theorem}\label{result} Let $\mathcal{A}$ be an intersecting
subfamily of $\mathcal{P}_{k,r,n}$. Then
\begin{enumerate}[(a)]
\item $|\mathcal{A}| \leq {k-1 \choose r-1}\frac{(n-1)!}{(n-r)!} =
{n-1 \choose r-1}\frac{(k-1)!}{(k-r)!}$, \label{conditiona}
\item the bound in (\ref{conditiona}) is attained if and only if
  $\mathcal{A}$ is a star of $\mathcal{P}_{k,r,n}$.
\end{enumerate}
\end{theorem}

\section{Proof of the result}

We will prove Theorem~\ref{result} by extending the arguments in
\cite{KL,LW} to our more general setting. Recall that we are assuming
$k \leq n$ and that Theorem~\ref{LM} settles our problem for the case
$r = k$, so we will only consider $r \leq k-1$. We will abbreviate
$\mathcal{P}_{k,r,n}$ to $\mathcal{P}$.

For convenience, we shall use `$\ourmod$' to represent the usual
\emph{modulo operation} with the exception that for any non-zero
integers $a$ and $b$ the value of $ba \ourmod a$ will be $a$, rather
than $0$.

Let $X$ be a set, and let $m = |X|$. A bijection $\sigma : X
\rightarrow [m]$ is called a \emph{cyclic ordering of $X$}; all the
elements in $X$ are arranged in a cycle, and $x \in X$ is the
$\sigma(x)$-th element in the cycle. If $\sigma$ is a cyclic
ordering of $X$ and the elements of a subset $A$ of $X$ are
numbered consecutively, in the cyclic sense, by $\sigma$, then we say
that $A$ \emph{meets} $\sigma$.

Katona's cycle method is based on the following fundamental result.

\begin{lemma}[\cite{K}]\label{cyclelemma} Let $X$ be a set
  of size at least $2r$, and let $\sigma$ be a cyclic ordering of
  $X$. Let $\mathcal{B} := \{B \in \binom{X}{r} \colon B \textrm{
    meets } \sigma \}$, and let $\mathcal{A}$ be an intersecting
  subfamily of $\mathcal{B}$. Then $|\mathcal{A}| \leq r$. Moreover,
  if $|X| > 2r$, then $|\mathcal{A}| = r$ if and only if $\mathcal{A}$
  is a star of $\mathcal{B}$.
\end{lemma}

The union of all sets in $\mathcal{P}$ is the Cartesian product $[k]
\times [n]$. We say that a cyclic ordering $\sigma$ of $[k] \times
[n]$ is \emph{$r$-good} if every set of $r$ elements $(x_1,y_1),
\dots, (x_r,y_r)$ of $[k]\times [n]$ that are numbered consecutively,
in the cyclic sense, by $\sigma$ are such that $x_1, \dots, x_r$ are
distinct and $y_1, \dots, y_r$ are distinct. In an $r$-good cyclic ordering any $r$ consecutive elements form a generalized permutation in
$\mathcal{P}$.

We will define a cyclic ordering of $[k] \times [n]$ that is $r$-good
for all $r \in [k-1]$. (It is interesting to note that no such cyclic
ordering exists if $r=k=n$.) Let $\tau \colon [k]\times [n]
\rightarrow [kn]$ be defined by
\begin{equation*} 
\tau(x,y) := k(y - x \ourmod n) + x.
\end{equation*}
As one can immediately see from the following example with $k = 5$ and
$n = 7$, where each element $(x,y)$ of $[k] \times [n]$ is given the
label $\tau(x,y)$ shown in bold superscript, $\tau$ is $(k-1)$-good,
and hence $\tau$ is $r$-good for all $r \in [k-1]$.

\begin{center}
\begin{tabular*}{0.75\textwidth}{@{\extracolsep{\fill} } lllll }
$(1,7)^{\bf 31}$ & $(2,7)^{\bf 27}$ & $(3,7)^{\bf 23}$ & $(4,7)^{\bf 19}$ &  $(5,7)^{\bf 15}$ \\
$(1,6)^{\bf 26}$ & $(2,6)^{\bf 22}$ &$(3,6)^{\bf 18}$ &$(4,6)^{\bf 14}$ & $(5,6)^{\bf 10}$\\
$(1,5)^{\bf 21}$ & $(2,5)^{\bf 17}$ &$(3,5)^{\bf 13}$ &$(4,5)^{\bf 9}$ & $(5,5)^{\bf 5}$ \\
$(1,4)^{\bf 16}$ & $(2,4)^{\bf 12}$ &$(3,4)^{\bf 8}$ &$(4,4)^{\bf 4}$ & $(5,4)^{\bf 35}$ \\
$(1,3)^{\bf 11}$ & $(2,3)^{\bf 7}$ &$(3,3)^{\bf 3}$ &$(4,3)^{\bf 34}$ & $(5,3)^{\bf 30}$ \\
$(1,2)^{\bf 6}$  & $(2,2)^{\bf 2}$ &$(3,2)^{\bf 33}$ &$(4,2)^{\bf 29}$ & $(5,2)^{\bf 25}$ \\
$(1,1)^{\bf 1}$  & $(2,1)^{\bf 32}$ &$(3,1)^{\bf 28}$ &$(4,1)^{\bf 24}$ & $(5,1)^{\bf 20}$ \\
\end{tabular*}
\end{center}

Let $S_n$ denote the set of all bijections from $[n]$ to $[n]$. For
any $(\phi, \psi) \in S_k \times S_n$, define $\tau_{\phi,\psi} : [k]
\times [n] \rightarrow [kn]$ by
\[
\tau_{\phi, \psi}(x,y) := \tau(\phi^{-1}(x),\psi^{-1}(y))
\]
(i.e.~$\tau_{\phi, \psi}(\phi(i),\psi(j)) := \tau(i,j)$). Note that
$\tau_{\phi,\psi}$ is a cyclic ordering of $[k] \times [n]$ and let
\[
T_{k,n} := \{\tau_{\phi, \psi} \colon (\phi, \psi) \in S_k \times S_n\}.
\]
Further, for any $(\phi, \psi) \in S_k \times S_n$,
define $f_{\phi, \psi} : [k] \times [n] \rightarrow [k] \times [n]$ by
\[
f_{\phi, \psi}(x,y) := (\phi(x), \psi(y)).
\]

\begin{lemma}\label{orderings} 
  For all $(\phi,\psi) \in S_k \times S_n$ the ordering
  $\tau_{\phi,\psi}$ is an $r$-good cyclic ordering of $[k] \times
  [n]$.
\end{lemma}
\textbf{Proof.} Suppose $\tau_{\phi,\psi}$ is not an $r$-good cyclic
ordering. Then there exist two distinct elements $(a_1,b_1)$ and
$(a_2,b_2)$ of $[k]\times [n]$ such that
\[
\tau_{\phi,\psi}(a_2,b_2) = ( \tau_{\phi,\psi}(a_1,b_1) + p ) \ourmod kn
\]
for some $p \in [r-1]$, with either $a_1 = a_2$ or $b_1 = b_2$. If $a_1 = a_2$, then
\[
\tau(\phi^{-1}(a_1), \psi^{-1}(b_2)) = \left( \tau(\phi^{-1}(a_1), \psi^{-1}(b_1)) + p \right) \ourmod kn,
\]
but this contradicts the definition of $\tau$. Similarly, we cannot have $b_1 = b_2$. \qed
\\

Let $Z$ be a set and $\sigma$ be a cyclic ordering of $Z$.  Let $m$ be
an integer with $2 \leq m \leq |Z|$ and suppose that $z_1, \dots, z_m$
are distinct elements of $Z$. If $\sigma(z_{i+1}) = \sigma(z_i) + 1
\ourmod |Z|$ for each $i \in [m-1]$, then we say that the tuple $(z_1,
\dots, z_m)$ is an \emph{$m$-interval of $\sigma$}, and we call
$\{z_1, \dots, z_m\}$ the \emph{set corresponding to $(z_1, \dots,
  z_m)$}. If $1 \leq m_1 \leq m_2 \leq m$ and $\ell = m_2 - m_1 + 1$, then
we call the $\ell$-interval $(z_{m_1}, \dots, z_{m_2})$ of $\sigma$ an
\emph{$\ell$-subinterval of $(z_1, \dots, z_m)$}.

\begin{lemma}\label{count}
Each member of $\mathcal{P}$ meets exactly $r!(k-r)! (n-r)! kn$
members of $T_{k,n}$.  
\end{lemma}
\textbf{Proof.}  Let $P, Q \in \mathcal{P}$. Clearly, $Q =
\{f_{\pi,\rho}(x,y) \colon (x,y) \in P\}$ for some $(\pi,\rho) \in S_k
\times S_n$.

Let $\tau_{\phi, \psi} \in T_{k,n}$ such that $P$ meets $\tau_{\phi,
  \psi}$. Then $Q$ meets $\tau_{\phi, \psi} \circ ({f_{\pi^{-1}, \rho^{-1}}})$. For any $(x,y) \in S_k \times S_n$
\begin{align*}
\tau_{\phi, \psi} \circ ({f_{\pi^{-1}, \rho^{-1}}})(x,y) &= \tau_{\phi, \psi}(\pi^{-1}(x),\rho^{-1}(y)) \\
                                                &= \tau(\phi^{-1} \circ \pi^{-1}(x),\psi^{-1} \circ \rho^{-1}(y)) \\
                                                &= \tau((\pi \circ \phi)^{-1}(x),(\rho \circ \psi)^{-1}(y)). 
\end{align*}
Thus, since $(\pi \circ \phi)^{-1} \in S_k$ and $(\rho \circ \psi)^{-1} \in S_n$, we have
\[
\tau_{\phi, \psi} \circ ({f_{\pi^{-1}, \rho^{-1}}}) = \tau_{(\pi \circ \phi)^{-1}, (\rho \circ \psi)^{-1}} \in T_{k,n}.
\]
So $Q$ meets at least as many members of $T_{k,n}$ as $P$
does. Conversely, we can do this for every ordering that $Q$ meets,
thus $P$ and $Q$ meet the same number of members of $T_{k,n}$.

Each of the $k!n!$ members of $T_{k,n}$ contains exactly $kn$
$r$-intervals, and, by Lemma~\ref{orderings}, the sets corresponding to these $r$-intervals are members of $\mathcal{P}$. Thus, for each $\tau_{\phi, \psi} \in T_{k,n}$, the number of members of $\mathcal{P}$ that meet $\tau_{\phi, \psi}$ is $kn$. Since $|\mathcal{P}| = \binom{k}{r}\frac{n!}{(n-r)!}$, each member of $\mathcal{P}$ meets exactly
\[
\frac{k!n!kn}{\binom{k}{r}\frac{n!}{(n-r)!}} = r! (k-r)! (n-r)! kn
\]
members of $T_{k,n}$.\qed
\\

For each $\tau_{\phi, \psi} \in T_{k,n}$ the \emph{characteristic vector of $\tau_{\phi, \psi}$} is the length-$\left( \binom{k}{r}\frac{(n)!}{(n-r)!} \right)$ vector in which each position
corresponds to a member $P$ of $\mathcal{P}$, and the entry is $1$ if
$P$ meets $\tau_{\phi, \psi}$, and $0$ otherwise. Similarly, for any
$\mathcal{A} \subseteq \mathcal{P}$, the \emph{characteristic vector $\chi_\mathcal{A}$ of $\mathcal{A}$} is the length-$\left(\binom{k}{r}\frac{(n)!}{(n-r)!}\right)$ vector in which each position
corresponds to a member $P$ of $\mathcal{P}$, and the entry is $1$ if
$P \in \mathcal{A}$, and $0$ otherwise. We now have the tools to prove
Theorem~\ref{result}.
\\
\\
\textbf{Proof of Theorem~\ref{result}.}  Let $\mathcal{A}$ be an
intersecting subfamily of $\mathcal{P}$ of maximum size. Define a
matrix $M$ in which the rows are indexed by the members of
$\mathcal{P}$, the columns are indexed by the members $\tau_{\phi,
  \psi}$ of $T_{k,n}$, and the column for $\tau_{\phi, \psi}$ is the
characteristic vector of $\tau_{\phi, \psi}$. For any $\ell \in
\mathbb{N}$, let $\one_\ell$ denote the all ones vector of length $\ell$. By Lemma~\ref{count},
\[ M {\one_{|T_{k,n}|}}^T = r! (k-r)! (n-r)!kn {\one_{|\mathcal{P}|}}^T. \]

Define $\mathcal{A}_{\phi,\psi}$ to be the set of all the members of
$\mathcal{A}$ that meet $\tau_{\phi,\psi}$.  Then the
$\tau_{\phi,\psi}$-entry of $\chi_\mathcal{A} M$ is equal to
$|\mathcal{A}_{\phi,\psi}|$; by Lemma~\ref{cyclelemma}, this value is
no more than $r$. So
\begin{align}
(k!)(n!)r \geq \chi_{\mathcal{A}} M {\one_{|T_{k,n}|}}^T 
   = r! (k-r)! (n-r)! kn \chi_{\mathcal{A}}  {\one_{ |\mathcal{P}|}}^T 
   = r! (k-r)! (n-r)! kn |\mathcal{A}|, \label{inequality}
\end{align}
which implies that
\[
|\mathcal{A}| \leq \frac{(k!)(n!)r}{r! (k-r)! (n-r)!kn} = \binom{k-1}{r-1}\frac{ (n-1)!}{(n-r)!}.
\]
This gives the first statement of Theorem~\ref{result}. 

The intersecting family $\{P \in \mathcal{P} \colon (1,1) \in P\}$
meets this bound, so the size of $\mathcal{A}$ is
$\binom{k-1}{r-1}\frac{ (n-1)!}{(n-r)!}$. Thus, equality holds in
(\ref{inequality}), and $|\mathcal{A}_{\phi,\psi}| = r$ for each
$\tau_{\phi,\psi} \in T_{k,n}$.  So Lemma~\ref{cyclelemma} tells us
that for each $\tau_{\phi,\psi} \in T_{k,n}$ the set
$\mathcal{A}_{\phi,\psi}$ consists of those $r$ sets that meet
$\tau_{\phi,\psi}$ and contain a fixed element
$(x_{\phi,\psi},y_{\phi,\psi})$. Thus, for each $\tau_{\phi,\psi} \in
T_{k,n}$,
\begin{align}
\mathcal{A}_{\phi,\psi} = \{A \colon A \mbox{ corresponds to an $r$-subinterval of } L_{\phi,\psi}\}, \label{ref4}
\end{align}
where $L_{\phi,\psi}$ is the $(2r-1)$-interval of $\tau_{\phi,\psi}$ with middle entry $(x_{\phi,\psi},y_{\phi,\psi})$.

Let $\beta$ be the identity function from $[k]$ to $[k]$, and let
$\gamma$ be the identity function from $[n]$ to $[n]$. So $\tau =
\tau_{\beta,\gamma}$. We may assume that
$(x_{\beta,\gamma},y_{\beta,\gamma}) = (k,k)$. So
$\mathcal{A}_{\beta,\gamma}$ consists of the $r$ sets corresponding to
all the $r$-subintervals of the $(2r-1)$-interval 
\[
L_{\beta,\gamma} = ((k-r+1,k-r+1), \dots, (k,k), (1,2), \dots, (r-1,r)). 
\]  
Define
\[
I := \{(i,i) \colon i \in [k-1]\}, \qquad \bar{I} := ([k] \times [n]) \backslash (I \cup \{(k,k)\}).
\] 
If $P \subseteq I$, then $P$ does not intersect the set
$\{(k,k),(1,2), \dots, (r-1,r)\} \in \mathcal{A}_{\beta,\gamma}$; similarly, if $P \subseteq \bar{I}$, then $P$ does not intersect the set $\{(k-r+1,k-r+1), \dots, (k,k)\} \in \mathcal{A}_{\beta,\gamma}$. Thus, for each $A \in \mathcal{A}$, it is the case that
$A \nsubseteq I$ and $A \nsubseteq \bar{I}$, so 
\begin{eqnarray}
 1 \leq |A \cap I| \leq r-1, \qquad  1 \leq |A \cap \bar{I}| \leq r-1. \label{ref5}
\end{eqnarray}

Define the sets
\[
T' := \{\tau_{\pi,\rho} \in T_{k,n} \colon \pi(k) = \rho(k) = k\}, \qquad
T^* := \{\tau_{\pi,\rho} \in T' \colon \pi(i) = \rho(i), i = 1, \dots, k\}.
\]
Note for each $\tau_{\pi,\rho} \in T^*$ that
\begin{eqnarray} 
\{(\pi(i),\rho(i)) \colon (i,i) \in I\} = I, \qquad 
\{(\pi(i),\rho(j)) \colon (i,j) \in \bar{I}\} = \bar{I}. \label{ref6}
\end{eqnarray}
If $(x_{\pi,\rho},y_{\pi,\rho}) \in I$, then, by (\ref{ref6}), $I$ has an $r$-subset $R$ that corresponds to an $r$-subinterval of $L_{\pi,\rho}$, and hence $R \in \mathcal{A}$ by (\ref{ref4}), but this contradicts the first inequality in (\ref{ref5}). Similarly,  $(x_{\pi,\rho},y_{\pi,\rho}) \in \bar{I}$ contradicts the second inequality in (\ref{ref5}). So $(x_{\pi,\rho},y_{\pi,\rho}) = (k,k)$ for each $\tau_{\pi,\rho} \in T^*$.


Now suppose $(x_{\pi,\rho},y_{\pi,\rho}) \neq (k,k)$ for some
$\tau_{\pi,\rho} \in T'$. Then $L_{\pi,\rho}$ has an $r$-subinterval
which does not have $(k,k)$ as one of its entries. Let $B$ be the
set corresponding to this interval; according to (\ref{ref4}), $B \in \mathcal{A}$. By (\ref{ref5}), $1 \leq s := |B \cap I| \leq r-1$. Let $(a_1,a_1), \dots, (a_s,a_s)$ be the $s$ distinct elements of $B \cap I$. Define $a_{s+1}, \dots, a_k$ to be the $k-s$ distinct elements of $[k] \backslash \{a_1, \dots, a_s\}$. Since $(k,k) \notin B \cap I$, we may assume that $a_k = k$. 

Choose $(\pi^*,\rho^*) \in S_k \times S_n$
such that $\pi^*(i) = \rho^*(i) = a_i$ for each $i \in [k]$. So
$\tau_{\pi^*,\rho^*} \in T^*$ and hence
$(x_{\pi^*,\rho^*},y_{\pi^*,\rho^*}) = (k,k) = (a_k,a_k)$ (as shown
above). Therefore,
\[
L_{\pi^*,\rho^*} = \left((a_{k-r+1},a_{k-r+1}), \dots, (a_k,a_k),(a_1,a_2), \dots, (a_{r-1},a_r)\right),
\] 
and the $r$-set 
\[
C := \{(a_{k-r+s},a_{k-r+s}), \dots, (a_k,a_k), (a_1,a_2), \dots, (a_{s-1},a_s)\}
\] 
corresponds to an $r$-subinterval of $L_{\pi^*,\rho^*}$; by
(\ref{ref4}), $C \in \mathcal{A}$. Since $k - r + s > s$, the pairs $(a_{k-r+s},a_{k-r+s}), \dots, (a_{k-1},a_{k-1}), (a_k,a_k)$ are not in $B$. Further, $(a_i,a_{i+1}) \notin B$ for each $i \in [s-1]$ since $(a_i,a_i) \in B$. Thus $B$ and $C$ are not intersecting, but this is a contradiction since $B, C \in \mathcal{A}$. We conclude that
\begin{eqnarray} 
(x_{\pi,\rho},y_{\pi,\rho}) = (k,k) \textrm{ for every } \tau_{\pi,\rho} \in T'. \label{ref7}
\end{eqnarray}

Finally, let $A$ be a set $\{(x_1,y_1), \dots, (x_r,y_r)\}$ in
$\mathcal{P}$ that contains $(k,k)$. We may assume that $(x_r,y_r) =
(k,k)$. Let $(\pi,\rho) \in S_k \times S_n$ be such that $\pi(i+k-r) =
x_i$ and $\rho(i+k-r) = y_i$ for each $i \in [r]$. Then
$\tau_{\pi,\rho} \in T'$ and $A$ meets $\tau_{\pi,\rho}$. By
(\ref{ref7}) and (\ref{ref4}), $A \in \mathcal{A}$. Hence the
result. \qed


\newpage
\begin{thebibliography}{}
\bibitem{Borg} P. Borg, Intersecting and cross-intersecting
families of labeled sets, Electron. J. Combin. 15 (2008) N9.

\bibitem{Borg7} P. Borg, Intersecting families of sets and
permutations: a survey, in: Advances in Mathematics Research
(A.R. Baswell Ed.), Volume 16, Nova Science Publishers, Inc.,
2011, pp 283--299, available at http://arxiv.org/abs/1106.6144.

\bibitem{Borg1} P. Borg, Intersecting systems of signed sets,
Electron. J. Combin. 14 (2007) R41.

\bibitem{Borg2} P. Borg, On $t$-intersecting families of signed
sets and permutations, Discrete Math. 309 (2009) 3310--3317.

\bibitem{BH} F. Brunk and S. Huczynska, Some Erd\H{o}s-Ko-Rado theorems for injections, European J. Combin. 31 (2010) 839--860.

\bibitem{CK} P.J. Cameron and C.Y. Ku, Intersecting families of
permutations, European J. Combin. 24 (2003) 881--890.

\bibitem{D} D.E. Daykin, Erd\H os-Ko-Rado from Kruskal-Katona,
J. Combin. Theory Ser. A 17 (1974) 254--255.

\bibitem{DF1} M. Deza and P. Frankl, On the maximum number of
permutations with given maximal or minimal distance, J. Combin.
Theory Ser. A 22 (1977) 352--360.

\bibitem{DF} M. Deza and P. Frankl, The Erd\H os-Ko-Rado theorem -- 22
years later, SIAM J. Algebraic Discrete Methods 4 (1983) 419--431.

\bibitem{EKR} P. Erd\H os, C. Ko and R. Rado, Intersection
theorems for systems of finite sets, Quart. J. Math. Oxford (2)
12 (1961) 313--320.

\bibitem{F} P. Frankl, The shifting technique in extremal set
theory, in: C. Whitehead (Ed.), Combinatorial Surveys, Cambridge
Univ. Press, London/New York, 1987, pp. 81--110.

\bibitem{GM} C. Godsil and K. Meagher, A new proof of the
  Erd\H{o}s-Ko-Rado theorem for intersecting families of permutations,
  European J. Combin. 30 (2009) 404--414.

\bibitem{HM} A.J.W. Hilton and E.C. Milner, Some intersection
theorems for systems of finite sets, Quart. J. Math. Oxford (2)
18 (1967) 369--384.

\bibitem{K} G.O.H. Katona, A simple proof of the Erd\H os-Chao
Ko-Rado theorem, J. Combin. Theory Ser. B 13 (1972) 183--184.

\bibitem{Kr} J.B. Kruskal, The number of simplices in a
complex, in: Mathematical Optimization Techniques, University of
California Press, Berkeley, California, 1963, pp. 251--278.

\bibitem{KL} C.Y. Ku and I. Leader, An Erd\H os-Ko-Rado theorem for
partial permutations, Discrete Math. 306 (2006) 74--86.

\bibitem{Li} Y.-S. Li, A Katona-type proof for intersecting families
of permutations, Int. J. Contemp. Math. Sciences 3 (2008) 1261--1268.

\bibitem{LM} B. Larose and C. Malvenuto, Stable sets of maximal
size in Kneser-type graphs, European J. Combin. 25 (2004) 657--673.

\bibitem{LW} Y.-S. Li and J. Wang, Erd\H os-Ko-Rado-type
theorems for colored sets, Electron. J. Combin. 14 (2007) R1.

\bibitem{WZ} J. Wang and S.J. Zhang, An Erd\H{o}s-Ko-Rado-type theorem
  in Coxeter groups, European J. Combin. 29 (2008) 1112--1115.
\end{thebibliography}
\end{document}